\documentstyle[leqno,12pt,twoside]{article}
\date{}

\begin{document}

\begin{titlepage}
\title{\bf The Asymptotic Behaviour of the Sum of Negative Eigenvalues of a Self-Adjoint Operator  Given in Semi-Axis}
\author {\bf \"Ozlem Bak\d{s}i }
\pagenumbering{arabic}\setcounter{page}{1} \pagestyle{myheadings}
\baselineskip=18pt \maketitle

\centerline{Department of Mathematics,}

\centerline{ Faculty of Arts and Science, Y\i ld\i z Technical University }
\centerline {(34210), Davutpa\d{s}a, \.Istanbul, Turkey}
\begin{center} e-mail: baksi@yildiz.edu.tr\end{center}

\noindent
\begin{abstract}
In this work, we find the asymptotic formulas for the sum of the negative eigenvalues smaller than
$\>-\varepsilon\>\quad(\varepsilon >0)\> $ of a self-adjoint operator $\>L\>$ which is defined by
the following differential expression
$$\ell(y)=-(p(x)y'(x))'-Q(x)y(x)$$
\noindent with the boundary condition

$$y(0)=0$$

\noindent in the space $\>L_{2}(0,\infty ;H)$.
\end{abstract}

\noindent {\small  {\bf AMS Subj. Classification:}} 34B24, 47A10


\end{titlepage}

%\vskip 0.5 cm

\section {Introduction}
Let $\>H\>$ be an infinite dimensional separable Hilbert space. Let us consider the operator $\>L\>$ in the Hilbert space $\>L_{2}(0,\infty ;H)\>$
defined by the differential equation
$$\ell(y)=-(p(x)y'(x))'-Q(x)y(x)\eqno(1)$$
\noindent and with the boundary condition $\>y(0)=0$.

\noindent Let us assume the scalar function $\>p(x)\>$ and the operator function $\>Q(x)\>$ satisfy the following conditions:

\noindent {\bf p1)} For every $\>x\in [0,\>\infty )$, there are positive constants $\>c_{1},\>c_{2}\>$ such that
$$c_{1}\leq p(x)\leq c_{2}.$$

\noindent {\bf p2)} The function $\>p(x)\>$ has continuous and bounded derivative.

\noindent {\bf p3)} The function $\>p(x)\>$ is not decreasing in the interval $\>[0,\> \infty )$.

\noindent {\bf Q1)} For every $\>x \in [0,\> \infty )\>$ the operator $\>Q(x):H\rightarrow H \>$ is self-adjoint, compact and positive.

\noindent {\bf Q2)} The operator $\>Q(x)\>$ is monotone decreasing.

\noindent {\bf Q3)} $\>Q(x)\>$ is a continuous operator function with respect to the norm in $\>B(H)\>$ and
$$\lim\limits _{x\rightarrow \infty }\|Q(x)\|=0.$$
$\>D(L)\>$ denotes the set of all functions $\>y(x) \in L_{2}(0,\> \infty; H )\>$ satisfying the following conditions:

\noindent {\bf y1)} $\>y(x)\>$ and $\>y'(x)\>$ are absolute continuous with respect to the norm in the space $\>H\>$ in every finite interval $\>[0,\>a]\>$

\noindent {\bf y2)} $\> l(y)=-(p(x)y'(x))'-Q(x)y(x)\in L_{2}(0,\> \infty; H )\>$

\noindent {\bf y3)} $\>y(0)=0\>$ and,
$$(Ly)(x)=-(p(x)y'(x))'-Q(x)y(x). \qquad \qquad $$
\noindent It is proved that the operator $\>L:D(L)\longrightarrow L_{2}(0,\> \infty;H )\>$ is self-adjoint, semi bounded-below and the negative part of the spectrum of the operator $\>L\>$ is discrete \cite{Se}.
\noindent Let $\>-\lambda_{1}\leq -\lambda_{2}\leq \cdots \leq -\lambda_{n}\leq \cdots \>$ be negative eigenvalues of the operator  $\>L\>$. In
this work we find an asymptotic formula for the sum
$$ \sum\limits_{-\lambda_{i}<-\varepsilon}\lambda_{i}\qquad (\varepsilon >0),\qquad $$
\noindent as $\>\varepsilon\rightarrow +0$.

\noindent In \cite{Ad} and \cite{Ba}, the asymptotic formulas for the sum of the negative eigenvalues of second order differential operator with
scalar coefficient are calculated. In \cite{Se}, \cite{Sk}, \cite{Ae}, \cite{Ma}, \cite{Ab} the asymptotic behaviour of the number of the negative eigenvalues are investigated.

\section{Some Inequalities For the Sum of the\\ Eigenvalues}

Let $\>\alpha_{1}(x)\geq \alpha_{2}(x)\geq \cdots \geq \alpha_{j}(x)\geq \cdots \>$ be the eigenvalues of the operator

\noindent $Q(x):H\longrightarrow H$. Since the operator function $\>Q(x)\>$ is monotone decreasing, the functions $\> \alpha_{1}(x),
\alpha_{2}(x),\cdots ,\alpha_{j}(x), \cdots \>$ are also monotone decreasing, \cite{Ae}.

\noindent Moreover, since
$$\alpha_{1}(x)=\sup_{\| f \|=1} (Q(x)f,f),$$
\noindent \cite{Co} and
$$\|Q(x)\|=\sup_{\| f \|=1}| (Q(x)f,f)|=\sup_{\|f \|=1} (Q(x)f,f),$$
\noindent  \cite{Ly} then $\>\alpha_{1}(x)=\|Q(x)\|$.

\noindent On the other hand, since $\>lim_{x\rightarrow\infty}\alpha_{1}(x)=0$, then the function $\>\alpha_{1}\>$ has a continuous
inverse function defined in the interval $\>(0,\alpha_{1}(0)]\>$. Let
\begin{eqnarray*}
\qquad\qquad \psi_{j}(\varepsilon)= \sup \{x\in[0,\infty);\>\alpha_{j}(x)\geq\varepsilon\}  \qquad (j=1,2,\cdots) \qquad\qquad\quad(2)
\end{eqnarray*}
\noindent and $\>\psi_{1}\>$ denote the inverse function of $\>\alpha_{1}$. We consider the following operators:

\noindent 1) Let $\>L^{0}\>$ and $\>L'\>$ be operators in the space $\>L_{2}(0,\> \psi_{1}(\varepsilon);H),\>$ which are formed by expression
$\>(1)\>$ and with the boundary conditions
\begin{eqnarray*}  y(0)&=& y(\psi_{1}(\varepsilon ))=0\\
                   y'(0)&=& y'(\psi_{1}(\varepsilon ))=0,
\end{eqnarray*}
\noindent respectively. Here, $\>\varepsilon \in (0,\alpha_{1}(0)]$.

\noindent 2) $\>L_{i}\>$ and $\>L_{i}'\>$ be operators in the space $\>L_{2}(x_{i-1},\>x_{i};H)\>$ which are formed by expression $\>(1)\>$ and
with the boundary conditions
\begin{eqnarray*} y(x_{i-1})&=& y(x_{i})=0\\
                       y'(x_{i-1})&=& y'(x_{i})=0,
\end{eqnarray*}
\noindent respectively.

\noindent 3)  $\>L_{i(1)}\>$ be operator in the space $\>L_{2}(x_{i-1},\>x_{i};H)\>$ which is formed by the differential equation
$$-p(x_{i})y''(x)-Q(x_{i})y(x) \qquad$$
\noindent and with the boundary conditions $\> y(x_{i-1})= y(x_{i})=0$.

\noindent 4) Let $\>L'_{i(1)}\>$  be operator in the space $\>L_{2}(x_{i-1},\>x_{i};H)\>$ which is formed by the differential equation
$$-p(x_{i-1})y''(x)-Q(x_{i-1})y(x) \qquad$$
\noindent and with boundary conditions $\> y'(x_{i-1})= y'(x_{i})=0$.

\noindent Let us divide the interval $\>[0,\> \psi_{1}(\varepsilon )]\>$ by the intervals at the length
$$ \delta ={{\psi_{1}(\varepsilon )}\over {[|\psi_{1}^{a}(\varepsilon )|]+1}} \qquad \eqno(3)$$
Here, $\> a \in (0,1) \>$ is a constant number and $\>\varepsilon \>$ is any positive number satisfying the inequality $\> \psi_{1}^{a } (\varepsilon )\geq2$. And also $\>[|\psi_{1}^{a}(\varepsilon )|]\>$ shows exact part of $\>\psi_{1}^{a}(\varepsilon )$.

\noindent Let the partition points of the interval  $\>[0,\>\psi_{1}(\varepsilon )]\>$ be
$$0=x_{0}<x_{1}<\cdots<x_{M}=\psi_{1}(\varepsilon).$$
\noindent Let $\>N(\lambda ),\>N^{0}(\lambda ),\>N'(\lambda ),\>n_{i}(\lambda )\>$ and $\>n_{i(1)}(\lambda )\>$ be numbers of eigenvalues smaller
than $\>-\lambda \quad (\lambda > 0)\>$ of the operators $\>L,\>L^{0},\>L',\>L_{i}\>$ and $\>L_{i(1)}\>$, respectively. Let us write
$ n_{i},\> n_{i(1)}\>$ instead of $ n_{i}(\varepsilon),\> n_{i(1)}(\varepsilon)\>$, respectively.

\noindent \c{S}eng\"{u}l \cite{Se} proved that the inequalities
$$ N^{0}(\varepsilon )\leq N(\varepsilon )\leq  N'(\varepsilon ) \eqno(4)$$
\noindent are satisfied, if $\>Q(x)\>$ satisfies the conditions $\>Q1),\>Q2),\>Q3)\>$ and $\>p(x)\>$ satisfies the conditions $\>p1),\>p3)$.

 \noindent We want to show that the inequalities
$$ N^{0}(\lambda )\leq N(\lambda )\leq  N'(\lambda ) \qquad (\forall \lambda \in [\varepsilon,\>\infty ))\eqno(5)$$
\noindent are satisfied.
\noindent Let $\>u_{1},u_{2},\cdots,u_{n},\cdots\>$ be orthonormal eigenvectors corresponding to the eigenvalues $\>-\lambda_{1}, -\lambda_{2},\cdots, -\lambda_{n}, \cdots $. Let us consider the following operators:

\begin{eqnarray*}
\qquad \qquad \qquad\quad\> S&=&L+\lambda I \qquad \qquad\qquad\qquad\qquad\qquad\qquad\qquad\qquad\>\>\>\>(6)\\
S^{0}&=&L^{0}+\lambda I ,\quad S'=L'+\lambda I \qquad\qquad\qquad \qquad\qquad\quad\>\>\> (7)
\end{eqnarray*}
\noindent Here $\>I\>$ in (6) is identity operator in the space $\>L_{2}(0,\infty;H)\>$; $\>I\>$ in (7) is identity operator in the space $\>L_{2}(0,\psi_{1}(\varepsilon);H)$. We have
$$-\lambda_{1}\leq -\lambda_{2}\leq\cdots \leq -\lambda_{N(\lambda)}<-\lambda, \qquad\lambda_{N(\lambda)+1}\geq-\lambda. \eqno(8) $$
\noindent Since the eigenvalues smaller than $\>\lambda\>$ are $\>\mu_{i}=\lambda_{i}+\lambda\quad (i=1,2,\cdots),\>$ from (8)
$$-\mu_{1}\leq-\mu_{2}\leq\cdots\leq-\mu_{N(\lambda)}<0, \quad -\mu_{N(\lambda)+1}\geq\>0 \eqno(9)$$
\noindent is obtained.By the similar way we can show that the number of negative eigenvalues of the operators $\>S^{0} \>\rm{and}\> S'\>$ are $\>N^{0}(\lambda) \>\rm{and}\> N'(\lambda)$, respectively. Let
\begin{eqnarray*}
 -\mu_{(1)1}\leq -\mu_{(1)2}\leq\cdots\leq-\mu_{(1)N^{0}(\lambda)},\quad
 -\mu_{(2)1}\leq -\mu_{(2)2}\leq\cdots\leq-\mu_{(2)N'(\lambda)}\quad(10)
\end{eqnarray*}

\noindent be negative eigenvalues of the operators $\>S^{0} \> \rm{and}\> S'\>$ respectively. Let the orthonormal eigenvectors corresponding these eigenvalues be $\varphi_{1},\varphi_{2},\cdots,\varphi_{N^{0}{(\lambda)}}\\
\rm{and}\quad \psi_{1},\psi_{2},\cdots,\psi_{N'{(\lambda)}}\>$ respectively.

\newtheorem{theorem}{Theorem}[section]
\newtheorem{lemma}[theorem]{Lemma}
\newtheorem{exam}[theorem]{Example}

\begin{lemma}
If the operator function $\>Q(x)\>$ satisfies the conditions \\
$\>Q1),\>Q2),\>Q3)\>$ and the function $\>p(x)\>$ satisfies the conditions $\>p1),\>p2)\>$ then
$$ N(\lambda)\geq N^{0}(\lambda)\quad \Big(\forall\lambda\in(0,\infty)\Big) \qquad\qquad\eqno(11)$$
\end{lemma}

\noindent{\bf Proof: } To obtain a contradiction, we suppose that
$$ N(\lambda)< N^{0}(\lambda).\qquad\qquad\qquad$$
\noindent Then, there is a non-zero linear combination
$$ \varphi=\sum\limits_{i=1}^{N^{0}(\lambda)}\beta_{i}\varphi_{i}\eqno(12)$$
\noindent of the functions $\>\varphi_{1},\varphi_{2},\cdots,\varphi_{N^{0}{(\lambda)}}\>$ such that
$$ \Big(u_{i},\>\varphi\Big)_{(0,\psi_{1}(\varepsilon))}=\int\limits_{0}^{\psi_{1}(\varepsilon)}\Big(u_{i}(x),\>\varphi (x)\Big)\>dx=0\quad(i=1,2,\cdots,N(\lambda)) $$

\noindent By using $\>(12)\>$
\begin{eqnarray*}
\Big(S^{0}\varphi,\>\varphi\Big)_{(0,\psi_{1}(\varepsilon))}&=&\bigg(S^{0}\Big(\sum\limits_{i=1}^{N_{1}(\lambda)}\beta_{i}\varphi_{i}\Big),\> \sum\limits_{i=1}^{N_{1}(\lambda)}\beta_{i}\varphi_{i}\bigg)_{(0,\psi_{1}(\varepsilon))}\\
\\
&=&\bigg(\sum\limits_{i=1}^{N_{1}(\lambda)}\beta_{i}\mu_{(1)i}\varphi_{i},\> \sum\limits_{i=1}^{N_{1}(\lambda)}\beta_{i}\varphi_{i}\bigg)_{(0,\psi_{1}(\varepsilon))}\\
\\
&>&\sum\limits_{i=1}^{N_{1}(\lambda)}\mu_{(1)i} \big|\beta_{i}\big|^{2}\>=\>\alpha<0 \qquad\qquad\qquad\qquad\quad(13)
\end{eqnarray*}

\noindent In the similar way as proved in Glazman \cite{Gl}  there exists a vector function $\>\tilde{\varphi}\>$ which has the following properties:

\noindent \textbf{$\tilde{\varphi}1)$} The vector function $\>\tilde{\varphi}=\tilde{\varphi}(x)\>$ has second second order continuous derivative respect to the norm in the space $\>H\>$ in the interval $\>[0,\psi_{1}(\varepsilon)]$.\\
\\
\noindent \textbf{$\tilde{\varphi}2)$} $\tilde{\varphi}(x)\>$ is equal to zero outside of the interval $\>[a,b]\subset\big(0,\psi_{1}(\varepsilon)\big)$.\\
\\
\noindent \textbf{$\tilde{\varphi}3)$}  $\Big|\Big(S^{0}\tilde{\varphi}\>,\>\tilde{\varphi}\Big)_{(0,\psi_{1}(\varepsilon))}-\>\Big(S^{0}\varphi\>,\> \varphi \Big)_{(0,\psi_{1}(\varepsilon))}\Big|\> < \> -\frac{\alpha}{2}$\\
\\
\noindent \textbf{$\tilde{\varphi}4)$} $\Big(u_{i}\> ,\> \tilde{\varphi}\Big)_{(0,\psi_{1}(\varepsilon))}=\>0 \quad \big(i=1,2,\cdots,N(\lambda)\big).$ \\

\noindent As it is known,
$$ \inf_{y\in D(S),\> \|y\|_{(0,\infty)}=1\atop{ y\perp u_{i}\>(i=1,2,\cdots,N(\lambda))}} \Big(Sy,y\Big)_{(0,\infty)}\>=\> \mu_{N(\lambda)+1}$$
\noindent Therefore
\begin{eqnarray*}
\bigg(S^{0}\Big(\frac{\tilde{\varphi}}{\parallel\tilde{\varphi}\parallel}\Big),\> \Big(\frac{\tilde{\varphi}}{\parallel\tilde{\varphi}\parallel}\Big)\bigg)_{(0,\psi_{1}(\varepsilon))}&=& \bigg(S\Big(\frac{\tilde{\varphi}}{\parallel\tilde{\varphi}\parallel}\Big),\> \Big(\frac{\tilde{\varphi}}{\parallel\tilde{\varphi}\parallel}\Big)\bigg)_{(0,\infty)}\\
\\
&\geq& \mu_{N(\lambda)+1}\>\geq 0.
\end{eqnarray*}
\noindent By the last inequality,
$$ \Big(S^{0}\varphi,\>\varphi\Big)_{(0,\psi_{1}(\varepsilon))}\geq 0\eqno(14)$$
\noindent is obtained. By $\>(13) \quad \rm{and} \>\> (14)\>$
$$ \Big(S^{0}\tilde{\varphi},\>\tilde{\varphi}\Big)_{(0,\psi_{1}(\varepsilon))}-\Big(S^{0}\varphi,\>\varphi\Big)_{(0,\psi_{1}(\varepsilon))} =\Big(S^{0}\tilde{\varphi},\>\tilde{\varphi}\Big)_{(0,\psi_{1}(\varepsilon))}-\>\alpha \geq -\alpha \eqno(15)$$
\noindent is found. On the other hand this result in $\>(15)\>$ contradicts with the property $\tilde{\varphi}$3) . Hence
$$ N(\lambda)\geq N^{0}(\lambda).$$
\begin{lemma}
If the operator function $\>Q(x)\>$ satisfies the conditions $\>Q1),\>Q2),\>Q3)\>$ and function $\>p(x)\>$ satisfies the conditions $\>p1),\> p2)\>$ then $\>N(\lambda)\>\leq\>N'(\lambda)\>$  for all $\>\lambda \in [\varepsilon,\> \infty)).$
\end{lemma}

\noindent {\bf Proof:} Suppose for contradiction that $\>N(\lambda)\> >\> N'(\lambda)$. Then, there is a non-zero linear combination

\newcommand{\ud}{\mathrm {d}}
\begin{displaymath}
\qquad\qquad\qquad u=\sum\limits_{i=1}^{N(\lambda)} d_{i} u_{i} \qquad\qquad\qquad\qquad\qquad\qquad\qquad\qquad\qquad\quad(16)
\end{displaymath}

\noindent of the vector functions $\>u_{1},\>u_{2},\cdots,\>u_{N(\lambda)}\>$ such that
$$ \Big(\psi_{i},\>u\Big)_{(0,\psi_{1}(\lambda))}=\int\limits_{0}^{\psi_{1}(\lambda)} \Big(\psi_{i}(x),\>u(x)\Big)\ud x=0\quad (i=1,2,\cdots,
N'(\lambda)) $$

\noindent By using $\>(16)$
\begin{eqnarray*}
\Big(Su,\>u\Big)_{(0,\infty)}&=&\bigg(S\Big(\sum\limits_{i=1}^{N(\lambda)}d_{i}u_{i}\Big)\>,\>\sum\limits_{i=1}^{N(\lambda)} d_{i}u_{i}\bigg)_{(0,\infty)}\\
\\
&=& \bigg(\sum\limits_{i=1}^{N(\lambda)}d_{i}\mu_{i}u_{i}\>,\>\sum\limits_{i=1}^{N(\lambda)} d_{i}u_{i}\bigg)_{(0,\infty)}\>=\>\sum\limits_{i=1}^{N(\lambda)} \mu_{i}\big|d_{i}\big|^{2}< \>0 \qquad(17)
\end{eqnarray*}
\noindent is obtained. We can write the equation $\>(17)\>$ as
$$ \Big(Su,\>u\Big)_{(0,\infty)}=\int\limits_{0}^{\psi_{1}(\varepsilon)}\Big(S(u(x)),\>u(x)\Big)\ud x\>+\> \int\limits_{\psi_{1}(\varepsilon)}^{\infty}\Big(S(u(x)),\>u(x)\Big)\ud x\> < 0 \eqno(18)$$
\noindent Since
$$\int\limits_{\psi_{1}(\varepsilon)}^{\infty}\Big(S(u(x)),\>u(x)\Big)\ud x\geq 0$$
\noindent then we have
$$\int\limits_{0}^{\psi_{1}(\varepsilon)}\Big(S(u(x)),\>u(x)\Big)\ud x < 0. \eqno(19)$$
\noindent If we consider the equality
$$ \Big(u,\psi_{i}\Big)_{(0,\psi_{1}(\varepsilon))}=\int\limits_{0}^{\psi_{1}(\varepsilon)} \Big(u(x), \psi_{i}(x)\Big)_{(0,\infty)}\ud x=0\quad (i=1,2,\cdots,N'(\lambda)) $$
\noindent from $\>(19)\>$
$$ \inf_{y\in D(S),\> \|y\|_{(0,\psi_{1}(\varepsilon))}=1\atop{ y\perp \psi_{i}\>(i=1,2,\cdots,N'(\lambda))}} \int\limits_{0}^{\psi_{1}(\varepsilon)}\Big(S(y(x)),y(x)\Big)\ud x< 0 \eqno(20)$$
\noindent is obtained. From $\>(20)$
$$ \inf_{y\in D(S),\> \|y\|_{(0,\psi_{1}(\varepsilon))}=1\atop{y'(0)=y'(\psi_{1}(\varepsilon))=0,\> y\perp \psi_{i}\>(i=1,2,\cdots,N'(\lambda))}} \int\limits_{0}^{\psi_{1}(\varepsilon)}\Big(S(y(x)),y(x)\Big)\ud x< 0 \eqno(21)$$
\noindent is found. By $\>(21)$
$$ \inf_{y\in D(S'),\> \|y\|_{(0,\psi_{1}(\varepsilon))}=1\atop{ y\perp \psi_{i}\>(i=1,2,\cdots,N'(\lambda))}} \int\limits_{0}^{\psi_{1}(\varepsilon)}\Big(S'(y(x)),y(x)\Big)\ud x< 0 \eqno(22)$$
\noindent is obtained. On the other hand, we have
$$ \inf_{y\in D(S'),\> \|y\|_{(0,\psi_{1}(\varepsilon))}=1\atop{ y\perp \psi_{i}\>(i=1,2,\cdots,N'(\lambda))}} \int\limits_{0}^{\psi_{1}(\varepsilon)}\Big(S'(y(x)),y(x)\Big)\ud x = \mu_{(2)(N'(\lambda)+1)}\geq 0 \eqno(23)$$
\noindent This result contradicts with (22). Therefore $\>N(\lambda)\leq N'(\lambda)\>$.  $\Box$

\noindent Let $\> -\mu_{i(1)1}\leq -\mu_{i(1)2}\leq -\mu_{i(1)3}\leq \cdots \>$ be eigenvalues of the operator $\>L_{i(1)}\>$
and let we have the following equalities
$$ a_{j}(x,t)=\alpha_{j}(x)-p(x)(\frac{\pi t}{\delta})^{2} \qquad (j=1,2,\cdots) \eqno(24)$$
$$  b_{j}(\varepsilon,x)=\frac{\delta}{\pi}\sqrt{\frac{\alpha_{j}(x)-\varepsilon}{p(x)}} \qquad (j=1,2,\cdots) \eqno(25)$$
$$  \beta_{j}(\varepsilon,x)=\int\limits_{0}^{b_{j}(\varepsilon,x)}a_{j}(x,t)dt \qquad (j=1,2,\cdots) \eqno(26)$$
$$ \varphi_{i,j}(\varepsilon)=\min\{x_{i+1},\>\psi_{j}(\varepsilon)\}\qquad (i=1,2,\cdots,M-1). \eqno(27)$$
\begin{theorem}
If the operator function $\>Q(x)\>$ and the scalar function $\>p(x)\>$ satisfy the conditions $\>Q1)-Q3)\>\rm{and}\> p1)-p3)$, then we have
$$\sum\limits_{m=1}^{n_{i(1)}}\mu_{i(1)m}>\frac{1}{\delta }\sum\limits_{{j}\atop \alpha_{j}(x_{i})>\varepsilon}
\int\limits_{x_{i}}^{\varphi_{i,j}(\varepsilon)} \beta_{j}(\varepsilon,x )dx-3\sum\limits_{{j}\atop \alpha_{j}(0)>\varepsilon} \alpha_{j}(0)$$

\noindent for small positive values of $\>\varepsilon $.

\end{theorem}
\noindent {\bf Proof :} Let us consider the operator $\>L_{i(1)}\>$ which is formed by the differential expression
$$ -p(x_{i})y''(x)-Q(x_{i})y(x)$$

\noindent with the boundary conditions $\>y(x_{i-1})=y(x_{i})=0$.

\noindent We wish to obtain the eigenvalues of the operator $\>L_{i(1)}$. In order to find the eigenvalues, we will solve the eigenvalues problem
\begin{eqnarray*}
\qquad\qquad\qquad\qquad\qquad -d u''&=&\lambda u\\
u(a)&=&u(b)=0 \qquad\qquad\qquad\qquad\qquad\qquad\quad(28)
\end{eqnarray*}

\noindent in the space $\>L_{2}(a,b)$. Here, $\>a=x_{i-1},\>b=x_{i}\>$ and $\> d=p(x_{i})$. Moreover, $\>\gamma\>$ is an eigenvalue of the operator $\>Q(b):H\longrightarrow H$. The eigenvalues of boundary-value problem (28) are in the form
$$ \lambda_{n}=d.\Big(\frac{n\pi}{b-a}\Big)^{2} , \quad (n\in {\bf N}).$$

\noindent So, the eigenvalues of the operator $\>L_{i(1)}\>$ are of the form
$$\lambda_{n}-\gamma \>=\>p(x_{i})\Big(\frac{n\pi}{x_{i}-x_{i-1}}\Big)^{2} -\gamma .$$

\noindent Since the eigenvalues of the operator $\>Q(x):H\longrightarrow H\>$ are $\> \alpha_{1}(x)\geq \alpha_{2}(x)\geq\cdots\geq \alpha_{j}(x)\geq \cdots\>$ then the eigenvalues of the operator $\>L_{i(1)}\>$ are
$$ p(x_{i})\biggl({{m\pi } \over {x_{i}-x_{i-1}}}\biggr )^{2}-\alpha_{j}(x_{i}) \quad(m=1,2,\cdots;\>j=1,2,\cdots),$$
\noindent therefore $\>n_{i(1)}\>$ is the number of pairs $\>(m,j)\qquad (m,j\geq1)\>$ satisfying the inequality
$$ p(x_{i})\bigl({{m\pi } \over {\delta}}\bigr )^{2}-\alpha_{j}(x_{i})<-\varepsilon \qquad (\delta=x_{i}-x_{i-1}). \eqno(29)$$

\noindent By using (24), (25) and (29), we obtain
\begin{eqnarray*}
\qquad\qquad\qquad\sum\limits_{m=1}^{n_{i(1)}}\mu_{i(1)m}&=&\sum\limits_{{j}\atop \alpha_{j}(x_{i})>\varepsilon} \sum\limits_{{m}\atop
a_{j}(x_{i},m)>\varepsilon}a_{j}(x_{i},m)\\
\\
&\geq & \sum\limits_{{j}\atop \alpha_{j}(x_{i})>\varepsilon} \sum\limits_{m=1}^{[|b_{j}(\varepsilon,x_{i})|]-1}a_{j}(x_{i},m)
\qquad\qquad\qquad\qquad(30)
\end{eqnarray*}

\noindent For the sum $\>\sum\limits_{m=1}^{[|b_{j}(\varepsilon,x_{i})|]-1}a_{j}(x_{i},m)\>$ in (30)
\begin{eqnarray*}
\qquad\qquad\sum\limits_{m=1}^{[|b_{j}(\varepsilon,x_{i})|]-1}a_{j}(x_{i},m)&\geq & \int\limits_{1}^{b_{j}
(\varepsilon,x_{i})-2}a_{j}(x_{i},t)\ud t\\
\\
&=&
\int\limits_{0}^{b_{j}(\varepsilon,x_{i})}a_{j}(x_{i},t)dt-\int\limits_{0}^{1}a_{j}(x_{i},t)\ud t\\
\\
&-&\int\limits_{b_{j}(\varepsilon,x_{i})-2}^{b_{j}(\varepsilon,x_{i})}a_{j}(x_{i},t)\ud t\\
\\
&>& \int\limits_{0}^{b_{j}(\varepsilon,x_{i})}a_{j}(x_{i},t)dt-3\alpha_{j}(x_{i})\\
\\
&=& \beta_{j}(\varepsilon,x_{i})-3\alpha_{j}(x_{i})\qquad\qquad\qquad\qquad\quad (31)
\end{eqnarray*}

\noindent is obtained. If we consider that the functions $\>\beta_{j}(\varepsilon,x) \quad (j=1,2,\cdots)\>$ are decreasing, from(27), (30) and (31)

\begin{eqnarray*}
\sum\limits_{m=1}^{n_{i(1)}}\mu_{i(1)m}&>& \frac{1}{\delta}\sum\limits_{{j}\atop \alpha_{j}(x_{i})>\varepsilon}
\int\limits_{x_{i}}^{x_{i+1}}\beta_{j}(\varepsilon,x_{i})dx-3\sum\limits_{{j}\atop \alpha_{j}(0)>\varepsilon}\alpha_{j}(0)   \\
\\
&\geq& \frac{1}{\delta}\sum\limits_{{j}\atop \alpha_{j}(x_{i})>\varepsilon}
\int\limits_{x_{i}}^{\varphi_{i,j}(\varepsilon)}\beta_{j}(\varepsilon,x)dx-3\sum\limits_{{j}\atop \alpha_{j}(0)>\varepsilon}\alpha_{j}(0)
\end{eqnarray*}

\noindent is obtained. $\Box$

\begin{theorem}
If the operator function $\>Q(x)\>$ and the scalar function $\>p(x)\>$ satisfy  the conditions $\>Q1)-Q3),p1)-p3)$, then we have
$$\sum\limits_{i=1}^{N(\varepsilon)}\lambda_{i}>{1\over {\delta
}}\sum\limits_{j=1}^{l_{\varepsilon}}\int\limits_{0}^{\psi_{j}(\varepsilon)}\beta_{j}(\>\varepsilon,x
)dx-const.\sum\limits_{j=1}^{l_{\varepsilon}}\int\limits_{0}^{\delta
}\alpha_{j}^{\frac{3}{2}}(x)dx-const.\psi_{1}^{a}(\varepsilon)\sum\limits_{j=1}^{l_{\varepsilon}}\alpha_{j}(0 )$$
for small positive values of $\>\varepsilon$.

\noindent Here, $\>l_{\varepsilon}=\sum\limits_{\alpha_{j}(0 )\geq \varepsilon}1 \>$.

\end{theorem}
\noindent {\bf Proof :} We can easily show that $\>L_{i}< L_{i(1)}\>$. In the case, it is known that
$$n_{i}(\lambda)\geq n_{i(1)}(\lambda) \eqno(32)$$
\noindent \cite{Sm}. On the other hand, from variation principles of R. Courant \cite{Ch}, we have
$$ N^{0}(\lambda)\geq \sum\limits_{i=1}^{M}n_{i}(\lambda ). \eqno(33)$$
\noindent From (32) and (33)
$$ N^{0}(\lambda)\geq \sum\limits_{i=1}^{M}n_{i(1)}(\lambda ) \qquad (\lambda\geq\varepsilon) \eqno(34)$$
\noindent is obtained. From (5) and (34)
$$ N(\lambda)\geq \sum\limits_{i=1}^{M}n_{i(1)}(\lambda ) \qquad (\forall\lambda\geq\varepsilon) \eqno(35)$$
\noindent is found. By using (35), we can show that the inequality
$$ \sum\limits_{i=1}^{N(\varepsilon)}\lambda_{i}\geq \sum\limits_{i=1}^{M}\sum\limits_{m=1}^{n_{i(1)}}\mu_{i(1)m} \eqno(36)$$
\noindent is satisfied. By the Theorem 2.1 and (36)
\begin{eqnarray*}
\sum\limits_{i=1}^{N(\varepsilon)}\lambda_{i}&\geq& \sum\limits_{i=1}^{M-1}\biggl\{\frac{1}{\delta}\sum\limits_{{j}\atop
\alpha_{j}(x)>\varepsilon}
\int\limits_{x_{i}}^{\varphi_{i,j}(\varepsilon)}\beta_{j}(\varepsilon,\>x)dx-3\sum\limits_{j=1}^{l_\varepsilon} \alpha_{j}(0)\biggr\}\\
\\
&=& \frac{1}{\delta}\sum\limits_{{j}\atop \alpha_{j}(x_{i})>\varepsilon}\sum\limits_{i}
\int\limits_{x_{i}}^{\varphi_{i,j}(\varepsilon)}\beta_{j}(\varepsilon,\>x)dx-3(M-1)\sum\limits_{j=1}^{l_\varepsilon} \alpha_{j}(0)
\qquad\qquad\quad(37)
\end{eqnarray*}

\noindent is obtained. Since the functions $\>\alpha_{j}(x) \quad (j=1,2,\cdots)\>$ are decreasing, then we have
$$ \sum\limits_{{j}\atop \alpha_{j}(x_{i})>\varepsilon}\sum\limits_{i}
\int\limits_{x_{i}}^{\varphi_{i,j}(\varepsilon)}\beta_{j}(\varepsilon,\>x)dx=\sum\limits_{{j}\atop
\alpha_{j}(x_{1})>\varepsilon}\sum\limits_{{i}\atop\alpha_{j}(x_{i})>\varepsilon}
\int\limits_{x_{i}}^{\varphi_{i,j}(\varepsilon)}\beta_{j}(\varepsilon,\>x)dx. \eqno (38)$$

\noindent From (37) and (38)
$$\sum\limits_{i=1}^{N(\varepsilon)}\lambda_{i}\geq \frac{1}{\delta}\sum\limits_{{j}\atop
\alpha_{j}(x_{1})>\varepsilon}\sum\limits_{{i}\atop\alpha_{j}(x_{i})>\varepsilon}
\int\limits_{x_{i}}^{\varphi_{i,j}(\varepsilon)}\beta_{j}(\varepsilon,\>x)dx-3M\sum\limits_{j=1}^{l_\varepsilon} \alpha_{j}(0)\eqno(39)$$

\noindent is obtained. By using (27) on the rigth-hand side of inequality (39)
\begin{eqnarray*}
\sum\limits_{i=1}^{N(\varepsilon)}\lambda_{i}&\geq&\frac{1}{\delta}\sum\limits_{\alpha_{j}(x_{1})>\varepsilon}
\Biggl[\int\limits_{x_{1}}^{x_{2}}\beta_{j}(\varepsilon,\>x)dx+\int\limits_{x_{2}}^{x_{3}}
\beta_{j}(\varepsilon,\>x)dx+\cdots+\int\limits_{x_{i_{0}}}^{\psi_{j}(\varepsilon)} \beta_{j}(\varepsilon,\>x)dx\Biggr]\\
\\
&-&3M\sum\limits_{j=1}^{l_\varepsilon} \alpha_{j}(0)\qquad\qquad\qquad\qquad\qquad\qquad\qquad\qquad\qquad\qquad\>(40)
\end{eqnarray*}

\noindent is found. Here, $\> i_{0}\>$ is a natural number satisfying the following condition:
$$ x_{i_{0}}<\psi_{j}(\varepsilon)\leq x_{i_{0}+1}.$$
\noindent By using (27) and (40)
\begin{eqnarray*}
\sum\limits_{i=1}^{N(\varepsilon)}\lambda_{i}&\geq& \frac{1}{\delta}\sum\limits_{\psi_{j}(\varepsilon)>x_{1}}
\int\limits_{x_{1}}^{\varphi_{j}(\varepsilon)} \beta_{j}(\varepsilon,\>x)dx-3M\sum\limits_{j=1}^{l_\varepsilon} \alpha_{j}(0)\\
\\
&=&\frac{1}{\delta}\sum\limits_{j=1}^{l_\varepsilon}\int\limits_{0}^{\psi_{j}(\varepsilon)}
\beta_{j}(\varepsilon,\>x)dx-\frac{1}{\delta}\sum\limits_{\psi_{j}(\varepsilon)<x_{1}} \int\limits_{0}^{\psi_{j}(\varepsilon)}
\beta_{j}(\varepsilon,\>x)dx \\
\\
&-&\frac{1}{\delta}\sum\limits_{\psi_{j}(\varepsilon)\geq x_{1}} \int\limits_{0}^{x_{1}}
\beta_{j}(\varepsilon,\>x)dx-3M\sum\limits_{j=1}^{l_\varepsilon} \alpha_{j}(0)\\
\\
&=&\frac{1}{\delta}\sum\limits_{j=1}^{l_(\varepsilon)} \int\limits_{0}^{\psi_{j}(\varepsilon)}
\beta_{j}(\varepsilon,\>x)dx-\frac{1}{\delta}\sum\limits_{j=1}^{l_(\varepsilon)} \int\limits_{0}^{\varphi_{0,j}(\varepsilon)}
\beta_{j}(\varepsilon,\>x)dx\\
\\
&-& 3M\sum\limits_{j=1}^{l_\varepsilon} \alpha_{j}(0) \qquad\qquad\qquad\qquad\qquad\qquad\qquad\qquad\qquad\qquad(41)
\end{eqnarray*}

\noindent is obtained.
\noindent From (24), (25) and (26)
\begin{eqnarray*}
{1\over \delta}\>\beta_{j}(\varepsilon,x)&=&{1\over \delta}\int\limits_{0}^{b_{j}(\varepsilon, x)}\>\Big[\alpha_{j}(x)-p(x)({{\pi
t}\over{\delta}})^{2}\Big]dt \\
\\
&=&{1\over \delta}\>\alpha_{j}(x)b_{j}(\varepsilon, x)-\pi^{2}\>{{p(x)}\over{3\delta^{3}}}\>b_{j}^{3}(\varepsilon, x) \\
\\
&=&{1\over \delta}\>b_{j}(\varepsilon, x)\>\Big[\alpha_{j}(x)-\pi^{2}\>{{p(x)}\over{3\delta^{2}}}\>b_{j}^{2}(\varepsilon, x)\Big] \\
\\
&=&{1\over \pi}\>{\sqrt{{\alpha_{j}(x)-\varepsilon}\over
{p(x)}}}\>\biggl[\alpha_{j}(x)-\pi^{2}\>{{p(x)}\over{3\delta^{2}}}.{{\delta^{2}\>\bigl(\alpha_{j}(x)-\varepsilon\bigr)}\over{\pi^{2}\>p(x)}}\bigg]\\
\\
&=&{1\over \pi}\>{\sqrt{{\alpha_{j}(x)-\varepsilon}\over {p(x)}}}\>\biggl[{2\over 3}\>\alpha_{j}(x)+\frac{\varepsilon}{3}\biggr]< const.
\alpha_{j}^{3\over 2}(x) \qquad\qquad\qquad(42)
\end{eqnarray*}

\noindent is found for the expression $\>{1\over \delta}\>\beta_{j}(\varepsilon,x)$. From (27) and (42),
$${1\over \delta}\sum\limits_{j=1}^{l_\varepsilon}\int\limits_{0}^{\varphi_{0,j}(\varepsilon)}\beta_{j}(\varepsilon,x)dx<const.
\sum\limits_{j=1}^{l_\varepsilon}\int\limits_{0}^{\delta}\alpha_{j}^{3\over 2}(x)dx \eqno(43)$$
\noindent is obtained. From (3), (41) and (43)
$$\sum\limits_{i=1}^{N(\varepsilon)}\lambda_{i}>{1\over \delta}\sum\limits_{j=1}^{l_\varepsilon}
\int\limits_{0}^{\psi_{j}(\varepsilon)}\beta_{j}(\varepsilon,x)dx-const.\sum\limits_{j=1}^{l_\varepsilon}
\int\limits_{0}^{\delta}\alpha_{j}^{3\over 2}(x)dx-const.\psi_{1}^{a}(\varepsilon)\sum\limits_{j=1}^{l_\varepsilon}\alpha_{j}(0)$$

\noindent is found.$\Box$

\par Let $\>-\mu_{i(1)1}'\leq-\mu_{i(1)2}'\leq-\mu_{i(1)3}'\leq\cdots\>$ be eigenvalues of the operator $\>L_{i(1)}'\>$ and
$\>n_{i(1)}'(\lambda)\>$ be number of the eigenvalues smaller than $\>-\lambda \>\>(\lambda >0)\>$  of the operator $\>L_{i(1)}'$.
 Moreover, we will simply write $\>n_{i(1)}'\>$ instead of $\>n_{i(1)}'(\varepsilon)$.

\begin{theorem}
If the operator function $\>Q(x)\>$ and the scalar function $\>p(x)\>$ satisfy the conditions $\>Q1)-Q3),\>p1)-p3)\>$ then the inequality
$$\sum\limits_{m=1}^{n_{i(1)}'}\mu_{i(1)m}'\leq {1\over \delta}\sum\limits_{{j}\atop \alpha_{j}(x_{i-1})>\varepsilon} \int\limits_{x_{i-2}}^
{x_{i-1}}\beta_{j}(\varepsilon, x)dx+\sum\limits_{j=1}^{l_{\varepsilon}}\alpha_{j}(0)\quad (i=2, 3, \cdots )$$

\noindent is satisfied for the small positive values of $\>\varepsilon$.
\end{theorem}

\noindent {\bf Proof:} The eigenvalues of the operator $\>L_{i(1)}'\>$ are in the form
$$ p(x_{i-1})\biggl [ {{(m-1)\pi}\over {(x_{i}-x_{i-1})}}\biggr ]^{2}-\alpha_{j}(x_{i-1})\qquad (m=1, 2, \cdots ; j=1, 2, \cdots ).$$
\noindent Therefore $\>n_{i(1)}'\>$ is the number of the pairs $\>(m, j )\quad (m, j \geq 1)\>$ satisfying the inequality
$$p(x_{i-1})\biggl [ {{(m-1)\pi}\over {(x_{i}-x_{i-1})}}\biggr ]^{2}-\alpha_{j}(x_{i-1})<-\varepsilon . \eqno(44)$$

\noindent From (24), (25), and (44)
\begin{eqnarray*}
\sum\limits_{m=1}^{n_{i(1)}'}\mu_{i(1)m}'&=&\sum\limits_{{j}\atop \alpha_{j}(x_{i-1})>\varepsilon}\sum\limits_{{m}\atop
a_{j}(x_{i-1},m-1)>\varepsilon}a_{j}(x_{i-1},m-1) \\
\\
&=&\sum\limits_{{j}\atop \alpha_{j}(x_{i-1})>\varepsilon}\sum\limits_{m=1}^{[\vert b_{j}(\varepsilon , x_{i-1})\vert]+1} a_{j}(x_{i-1},m-1)
\qquad\qquad\qquad\qquad\quad(45)
\end{eqnarray*}

\noindent is found. It is easy to see that
\begin{eqnarray*}
\sum\limits_{m=1}^{[\vert b_{j}(\varepsilon , x_{i-1})\vert]+1}a_{j}(x_{i-1},m-1) &\leq& \alpha_{j}(x_{i-1})+\int\limits_{0}^{b_{j}(\varepsilon , x_{i-1})}a_{j}(x_{i-1},t)
\ud t\\
\\
&=&\alpha_{j}(x_{i-1})+\beta_{j}(\varepsilon ,x_{i-1}).\qquad\qquad\qquad\quad(46)
\end{eqnarray*}

\noindent We consider that the functions $\>\beta_{j}(\varepsilon ,x) \quad (j=1, 2, \cdots )\>$ are monotone decreasing, by (45) and (46),
\begin{eqnarray*}
\sum\limits_{m=1}^{n_{i(1)}'}\mu_{i(1)m}'&\leq& \sum\limits_{j=1}^{l_{\varepsilon}}\alpha_{j}(0)+{1\over \delta}\sum\limits_{{j}\atop
\alpha_{j}(x_{i-1})>\varepsilon}\int\limits_{x_{i-2}}^ {x_{i-1}}\beta_{j}(\varepsilon, x_{i-1})dx\\
\\
&<&{1\over \delta}\sum\limits_{{j}\atop \alpha_{j}(x_{i-1})>\varepsilon} \int\limits_{x_{i-2}}^ {x_{i-1}}\beta_{j}(\varepsilon,
x)dx+\sum\limits_{j=1}^{l_{\varepsilon}}\alpha_{j}(0)\qquad (i=1, 2, \cdots )
\end{eqnarray*}

\noindent is obtained.$\>\Box $

\par Let $\>n_{i}'(\lambda)\>$ be number of the eigenvalues smaller than $\>-\lambda\quad (\lambda >0)\>$ of the operator $\>L_{i}'\>$,
$\>-\mu_{1}'\leq -\mu_{2}'\leq -\mu_{3}'\leq \cdots\>$ be eigenvalues of the operator $\>L_{1}'\>$ and  $\>n_{i}'(\varepsilon)=n_{i}'$.

\begin{theorem}
If the operator function $\>Q(x)\>$ and the scalar function $\>p(x)\>$ satisfy the conditions $\>Q1)-Q3),\> \rm{and}\> p1)-p3)$, then we have
$$\sum\limits_{i=1}^{N(\varepsilon)}\lambda_{i} < \sum\limits_{m=1}^{n_{1}'}\mu_{m}'+{{1}\over
{\delta}}\sum\limits_{j=1}^{l_{\varepsilon}}\int\limits_{0}^{\psi_{j}(\varepsilon)}\beta_{j}(\varepsilon , x) dx+{{\psi_{1}(\varepsilon)}\over
\delta}\sum\limits_{j=1}^{l_{\varepsilon}}\alpha_{j}(0)$$

\noindent for the small values of $\>\varepsilon\>$.
\end{theorem}

\noindent {\bf Proof :} We can easily show that $\>L_{i}'>L_{i(1)}'$. In this case we have
$$n_{i}'(\lambda)\leq n_{i(1)}'(\lambda) \eqno(47)$$
\noindent \cite{Sm}. On the other hand, from variation principles of R. Courant \cite{Ch}, we have
$$N'(\lambda)\leq \sum\limits_{i=1}^{M}n_{i}'(\lambda).\eqno (48)$$
\noindent From (47) and (48),
$$ N'(\lambda)\leq \sum\limits_{i=2}^{M}n_{i(1)}'(\lambda)+n_{1}'(\lambda)\eqno(49)$$
\noindent is obtained. From (5) and (49)
$$ N(\lambda)\leq \sum\limits_{i=2}^{M}n_{i(1)}'(\lambda)+n_{1}'(\lambda) \qquad (\forall \lambda \geq \varepsilon) \eqno(50)$$
\noindent is found. By using (50),we have
$$ \sum\limits_{i=1}^{N(\varepsilon)}\lambda_{i}\leq
\sum\limits_{i=2}^{M}\sum\limits_{m=1}^{n_{i(1)}'}\mu_{i(1)m}'+\sum\limits_{m=1}^{n_{1}'}\mu_{m}'.\eqno(51)$$
\noindent By using Theorem 2.3 and (51)
\begin{eqnarray*}
\sum\limits_{i=1}^{N(\varepsilon)}\lambda_{i}&\leq&\sum\limits_{m=1}^{n_{1}'}\mu_{m}'+
{1\over \delta}\sum\limits_{i=2}^{M}\sum\limits_{{j}\atop \alpha_{j}(x_{i-1})>\varepsilon}\int\limits_{x_{i-2}}^ {x_{i-1}}
\beta_{j}(\varepsilon, x)dx+ M\sum\limits_{j=1}^{l_{\varepsilon}}\alpha_{j}(0)\\
\\
&=&\sum\limits_{m=1}^{n_{1}'}\mu_{m}'+{1\over \delta}\sum\limits_{{j}\atop\alpha_{j}(x_{i-1})>\varepsilon}\sum\limits_{i\geq 2}\>
\int\limits_{x_{i-2}}^{x_{i-1}}\beta_{j}(\varepsilon , x) dx
+M\sum\limits_{j=1}^{l_{\varepsilon}}\alpha_{j}(0)\qquad\quad(52)
\end{eqnarray*}
\noindent is found. Since the functions $\>\alpha_{j}(x)\quad (j=1,2,\cdots )\>$ are monotone decreasing, then we have
\begin{eqnarray*}
\qquad\quad\sum\limits_{{j}\atop\alpha_{j}(x_{i-1})>\varepsilon}\sum\limits_{i\geq 2}
\int\limits_{x_{i-2}}^{x_{i-1}}\beta_{j}(\varepsilon , x)dx=\sum\limits_{{j}\atop\alpha_{j}(x_{1})>\varepsilon}
\sum\limits_{{i\geq 2}\atop\alpha_{j}(x_{i-1})>\varepsilon}
\int\limits_{x_{i-2}}^{x_{i-1}}\beta_{j}(\varepsilon , x)\ud x. \qquad\quad(53)
\end{eqnarray*}
\noindent From (52) and (53)
\begin{eqnarray*}
\sum\limits_{i=1}^{N(\varepsilon)}\lambda_{i}&<&
\sum\limits_{m=1}^{n_{1}'}\mu_{m}'+{1\over \delta}\sum\limits_{{j}\atop \alpha_{j}(x_{1})>\varepsilon }
\sum\limits_{{i\geq 2}\atop\alpha_{j}(x_{i-1})>\varepsilon}
\int\limits_{x_{i-2}}^{x_{i-1}}\beta_{j}(\varepsilon , x)dx+M\sum\limits_{j=1}^{l_{\varepsilon}}\alpha_{j}(0)\\
\\
&=&\sum\limits_{m=1}^{n_{1}'}\mu_{m}'+{1\over\delta}\sum\limits_{{j}\atop \alpha_{j}(x_{1})>\varepsilon }\Bigl[\int\limits_{0}^{x_{1}}
\beta_{j}(\varepsilon , x)dx+\int\limits_{x_{1}}^{x_{2}}\beta_{j}(\varepsilon , x)dx\\
\\
&+&\cdots+\int\limits_{x_{i_{0}-1}}^{x_{i_{0}}}
\beta_{j}(\varepsilon , x) dx \Bigr ]
+M\sum\limits_{j=1}^{l_{\varepsilon}}\alpha_{j}(0)
\end{eqnarray*}
\noindent is obtained. Here, $\>i_{0}\>$ is a natural number satisfying the conditions
$$ \alpha_{j}(x_{i_{0}})>\varepsilon, \qquad \alpha_{j}(x_{i_{0}+1})\leq \varepsilon. \eqno(54)$$
\noindent From (2)
$$ x_{i_{0}}\leq \psi_{j}(\varepsilon). \eqno(55)$$
\noindent From (54) and (55)
$$ \sum\limits_{i=1}^{N(\varepsilon)}\lambda_{i}<\sum\limits_{m=1}^{n_{1}'}\mu_{m}'+
{1\over \delta}\sum\limits_{j=1}^{l_{\varepsilon}}\int\limits_{0}^{\psi_{j}(\varepsilon)}\beta_{j}(\varepsilon , x) dx +
{{\psi_{1}(\varepsilon)}\over {\delta}}\sum\limits_{j=1}^{l_{\varepsilon}}\alpha_{j}(0)$$

\noindent is found.$\Box$

\noindent Let
$$\delta_{i}={{\delta_{i-1}}\over {[\vert \delta_{i-1}\psi_{1}^{(i+1)a-1}(\varepsilon)\vert]+1}}\> ,\qquad (i=1,2,\cdots ;
\delta_{0}=\delta)\eqno(56)$$
$$a_{j(i)}(x,t)=\alpha_{j}(x)-p(x)\biggl({{\pi t}\over {\delta_{i}}}\biggr)^{2}\> ,$$
$$b_{j(i)}(\varepsilon, x)={{\delta_{i}}\over \pi}\sqrt{{{\alpha_{j}(x)-\varepsilon}\over {p(x)}}}\> ,$$
$$\beta_{j(i)}(\varepsilon,\>x)=\int\limits_{0}^{b_{j(i)}(\varepsilon, x)}a_{j(i)}(x,t)\ud t\>,$$
$$\varphi_{j}(\delta_{i},\varepsilon)=\min\{\delta_{i},\psi_{j}(\varepsilon)\}\qquad (i=0,1,2,\cdots ). \eqno(57)$$
\noindent Let $\>L_{(i)}\>$ be operator in the space $\>L_{2}(0,\delta_{i};H)\>$ which is formed by the expression (1) and with the boundary
condition
$$ y'(0)=y'(\delta_{i})=0. \eqno(58)$$
\noindent Moreover, let $\>L_{(i)}^{(0)}\>$ be operator  which is formed by the expression
$$\>-p(0)y''(x)-Q(0)y(x)\>$$
\noindent and with the boundary condition (58).

\noindent Let $\>-\mu_{(i)1}\leq -\mu_{(i)2}\leq \cdots\>$ and $\>-\mu_{(i)1}^{(0)}\leq -\mu_{(i)2}^{(0)}\leq \cdots\>$
be eigenvalues smaller than $\>-\lambda,\quad (\lambda >0)\>$ of the operators $\>L_{(i)}\>$ and $\>L_{(i)}^{(0)}\>$, respectively.

\noindent Moreover, let $\>n_{(i)}(\lambda)\>$ and $\>n_{(i)}^{(0)}(\lambda)\>$ be numbers of the eigenvalues smaller than $\>-\lambda,\>
(\lambda >0)\>$ of the operators $\>L_{(i)}\>$ and $\>L_{(i)}^{(0)}\>$, respectively.

Since $\>L_{(i)}\geq L_{(i)}^{(0)}\>$, then we have
$$ n_{(i)}(\lambda)\leq n_{(i)}^{(0)}(\lambda), \eqno(59)$$
\noindent \cite{Sm}. By using (59), we can show that
$$\sum\limits_{m=1}^{n_{(i)}}\mu_{(i)m}\leq\sum\limits_{m=1}^{n_{(i)}^{(0)}}\mu_{(i)m}^{(0)}.\eqno(60)$$
\noindent Here, $\>n_{(i)}=n_{(i)}(\varepsilon),\quad n_{(i)}^{(0)}=n_{(i)}^{(0)}(\varepsilon)$. $\>\delta_{-1}=\psi_{1}(\varepsilon)\>$ and from
the formula (56)
\begin{eqnarray*}{{\delta_{i-1}}\over {\delta_{i}}}&=&[\vert \delta_{i-1}\psi_{1}^{(i+1)a-1}(\varepsilon)\vert]+1\leq
\delta_{i-1}\psi_{1}^{(i+1)a-1}(\varepsilon)+1\\
\\
&=&{{\delta_{i-2}}\over {[\vert \delta_{i-2}\psi_{1}^{ia-1}(\varepsilon)\vert]+1}}\psi_{1}^{(i+1)a-1}(\varepsilon)+1\\
\\
&<&{{\delta_{i-2}}\over { \delta_{i-2}\psi_{1}^{ia-1}(\varepsilon)}}\psi_{1}^{(i+1)a-1}(\varepsilon)+1\\
\\
&=&\psi_{1}^{a}(\varepsilon)+1\qquad (i=1,2,\cdots)
\end{eqnarray*}

\noindent is obtained. From the last relation, we find
$${{\delta_{i-1}}\over {\delta_{i}}}<2\psi_{1}^{a}(\varepsilon), \qquad (i=1,2,\cdots )\eqno(61)$$
\noindent for the values of $\>\varepsilon\>$ satisfying the inequality $\>\psi_{1}^{a}(\varepsilon)>2$.

\begin{theorem}
If the operator function $\>Q(x)\>$ and the scalar function $\>p(x)\>$ satisfy the conditions $\>Q1)-Q3),\> \rm{and}\> p1)-p3)$, then we have
$$\sum\limits_{i=1}^{N(\varepsilon)}\lambda_{i}<{1\over
\delta}\sum\limits_{j=1}^{l_{\varepsilon}}\int\limits_{0}^{\psi_{j}(\varepsilon)}\beta_{j}(\varepsilon,x)dx+\mbox{const.} \sum\limits_{j=1}^{l_{\varepsilon}}\int\limits_{0}^{\delta}
\alpha_{j}^{{3\over 2}}(x)dx+\mbox{const.}\psi_{1}^{a}(\varepsilon)\sum\limits_{j=1}^{l_{\varepsilon}}\alpha_{j}(0)$$
\noindent for small positive values of $\>\varepsilon$.
\end{theorem}
\noindent {\bf Proof :} By the similar way to the proof of Theorem 2.6, the following inequality
\begin{eqnarray*}
\sum\limits_{m=1}^{n_{1}'}\mu_{m}'&<&\sum\limits_{m=1}^{n_{1}}\mu_{(1)m}+{1\over
{\delta_{1}}}\sum\limits_{\psi_{j}(\varepsilon)<\delta_{0}}\int\limits_{0}^{\psi_{j}(\varepsilon)}\beta_{j(1)}(\varepsilon ,x)\ud x \nonumber \\
\\
&+&{1\over {\delta_{1}}}\sum\limits_{\psi_{j}(\varepsilon)\geq \delta_{0}}\int\limits_{0}^{\delta_{0}}\beta_{j(1)}(\varepsilon ,x)\ud x+
{{\delta_{0}}\over {\delta_{1}}}\sum\limits_{j=1}^{l_{\varepsilon}}\alpha_{j}(0)\qquad\qquad\qquad\qquad\quad(62)
\end{eqnarray*}
\noindent can be proved. If we replace the equation (57) in (62), then we have
$$\sum\limits_{m=1}^{n_{1}'}\mu_{m}'<\sum\limits_{m=1}^{n_{1}}\mu_{(1)m}+{1\over {\delta_{1}}}
\sum\limits_{j=1}^{l_{\varepsilon}}
\int\limits_{0}^{\varphi_{j}(\delta_{0},\varepsilon)}\beta_{j(1)}(\varepsilon ,x)dx
+{{\delta_{0}}\over {\delta_{1}}}\sum\limits_{j=1}^{l_{\varepsilon}}\alpha_{j}(0).\eqno(63)$$
\noindent If we apply the inequality (63) for the eigenvalues of the operator $\>L_{(i)}\>$, then
$$\sum\limits_{m=1}^{n_{(i)}}\mu_{(i)m}<\sum\limits_{m=1}^{n_{(i+1)}}\mu_{(i+1)m}+{1\over {\delta_{i+1}}}
\sum\limits_{j=1}^{l_{\varepsilon}}
\int\limits_{0}^{\varphi_{j}(\delta_{i},\varepsilon)}\beta_{j(i+1)}(\varepsilon ,x)\ud x
+{{\delta_{i}}\over {\delta_{i+1}}}\sum\limits_{j=1}^{l_{\varepsilon}}\alpha_{j}(0)\eqno(64)$$
\noindent is obtained. From (61) and (64)
$$\sum\limits_{m=1}^{n_{(i)}}\mu_{(i)m}<\sum\limits_{m=1}^{n_{(i+1)}}\mu_{(i+1)m}+{1\over {\delta_{i+1}}}
\sum\limits_{j=1}^{l_{\varepsilon}}
\int\limits_{0}^{\varphi_{j}(\delta_{i},\varepsilon)}\beta_{j(i+1)}(\varepsilon ,x)dx
+2\psi_{1}^{a}(\varepsilon)\sum\limits_{j=1}^{l_{\varepsilon}}\alpha_{j}(0)\eqno(65)$$

\noindent is found. By using (45) and (46)
$$\sum\limits_{m=1}^{n_{(i+1)}^{(0)}}\mu_{(i+1)m}^{(0)}\leq \sum\limits_{j=1}^{l_{\varepsilon}}\Bigl(\alpha_{j}(0)+\beta_{j(i+1)}(\varepsilon ,
0\bigr)\Bigr)\eqno(66)$$
\noindent is obtained. Moreover, if we use the equation (42), then we get
$$\beta_{j(i+1)}(\varepsilon , x)\leq \mbox{const.}\delta_{i+1}\alpha_{j}^{3\over 2}(x).\eqno(67)$$
\noindent From (60), (66) and (67),
$$\sum\limits_{m=1}^{n_{(i+1)}}\mu_{(i+1)m}\leq \sum\limits_{j=1}^{l_{\varepsilon}}
\alpha_{j}(0)+\mbox{const.}\delta_{i+1}\sum\limits_{j=1}^{l_{\varepsilon}}\alpha_{j}^{3\over 2}(0)\eqno(68)$$
\noindent is obtained. By using inequality (56), we find
$$\delta_{i_{0}+1}\leq 1. \eqno(69)$$
\noindent Here, $\>i_{0} \in N\>$ is a constant satisfying the condition
$$i_{0}\geq {1\over a }-2$$
\noindent From (68) and (69), we get
$$\sum\limits_{m=1}^{n_{(i_{0}+1)}}\mu_{(i_{0}+1)m}\leq \mbox{const.}\sum\limits_{j=1}^{l_{\varepsilon}}\alpha_{j}(0).\qquad\qquad\eqno(70)$$

\noindent From (61), (63), (65) and (70),
\begin{eqnarray*}
\sum\limits_{m=1}^{n_{1}'}\mu_{m}'&\leq& \mbox{const.}\sum\limits_{j=1}^{l_{\varepsilon}}\alpha_{j}(0)+\sum\limits_{i=0}^{i_{0}}
{1\over {\delta_{i+1}}}\int\limits_{0}^{\varphi_{j}(\delta_{i},\varepsilon)}\beta_{j(i+1)}(\varepsilon
,x)\ud x\\
\\
&+& 2(i_{0}+1)\psi_{1}^{a}(\varepsilon)\sum\limits_{j=1}^{l_{\varepsilon}}\alpha_{j}(0) \qquad\qquad\qquad\qquad\qquad\qquad\qquad\qquad(71)
\end{eqnarray*}

\noindent is found. From (57), (67) and (71),
$$\sum\limits_{m=1}^{n_{1}'}\mu_{m}'< \mbox{const.}\sum\limits_{j=1}^{l_{\varepsilon}}\int\limits_{0}^{\delta}\alpha_{j}^{{3\over
2}}(x)dx+\mbox{const.}\psi_{1}^{a}(\varepsilon)\sum\limits_{j=1}^{l_{\varepsilon}}\alpha_{j}(0)\eqno(72)$$
\noindent is obtained. By the Theorem 2.6 and (72), we have
$$\sum\limits_{i=1}^{N(\varepsilon)}\lambda_{i}<{1\over
{\delta}}\sum\limits_{j=1}^{l_{\varepsilon}}\int\limits_{0}^{\psi_{j}(\varepsilon)}\beta_{j}(\varepsilon ,x)dx
+\mbox{const.}\sum\limits_{j=1}^{l_{\varepsilon}}\int\limits_{0}^{\delta}\alpha_{j}^{3\over
2}(x)dx+\mbox{const.}\psi_{1}^{a}(\varepsilon)\sum\limits_{j=1}^{l_{\varepsilon}}\alpha_{j}(0). $$

\noindent is obtained. $\Box$

\section{Asymptotic Formulas For The Sum Of \\ Negative Eigenvalues}

In this section, we find asymptotic formulas for the sum $\>\sum\limits_{-\lambda_{i}<-\varepsilon}\lambda_{i}\>$ as $\>\varepsilon\rightarrow
+0\>$.

\noindent Let us denote the functions of the form $\>\ln_{0}x=x, \quad \ln_{n}x=\ln(\ln_{n-1}x)\>$  by $\>\ln_{n}x\>\quad(n=0,1,2,\cdots )\>$ and
we suppose that the function $\>\alpha_{1}(x)=\Vert Q(x)\Vert \>$ satisfies the following condition:
\noindent \\
\\
\noindent {\bf{$\alpha 1)$}} There are a number $\>\xi>0\>$ and a natural number $\>n\geq1\>$ such that the function $\> \alpha_{1}(x)-(\ln_{n}x)^{-\xi}\>$
is neither negative nor monotone increasing in the interval $\>[b, \infty)\quad (b>0).\>$

\begin{theorem}
If the conditions $\>Q1)-Q3),\> p1)-p3)\> \rm{and}\> \alpha 1)\>$ are satisfied  and the series $\>\sum\limits_{j=1}^{\infty}[\alpha_{j}(0)]^{m}\>$ is convergent for a constant $\>m\in(0,\infty)\>$ , then the asymptotic formula
$$\sum\limits_{-\lambda_{i}<-\varepsilon}\lambda_{i}={1\over {3\pi}}\biggl[1 +O(e^{-\varepsilon^{-\beta}})\biggr]\sum\limits_{j}\int\limits_{\alpha_{j}(x)\geq
\varepsilon}\sqrt{{{\alpha_{j}(x)-\varepsilon}\over{p(x)}}}\>\Big(2\alpha_{j}(x)+\varepsilon\Big)\ud x$$
\noindent is satisfied as $\>\varepsilon\rightarrow +0\>$. Here, $\>\beta\>$ is a positive constant.
\end{theorem}

\noindent {\bf Proof:} By using Theorem 2.4 and Theorem 2.5, we have

$$\Bigl\vert \sum\limits_{i=1}^{N(\varepsilon)}\lambda_{i}-{1\over
{\delta}}\sum\limits_{j=1}^{l_{\varepsilon}}\int\limits_{0}^{\psi_{j}(\varepsilon)}\beta_{j}(\varepsilon ,x)dx \Bigr\vert <
\mbox{const.}l_{\varepsilon}(\delta +\psi_{1}^{a}(\varepsilon))$$

\noindent for the small positive values of $\>\varepsilon$. If we take $\> a={1\over 2}\>$ and consider (3)

$$\Bigl\vert \sum\limits_{i=1}^{N(\varepsilon)}\lambda_{i}-{1\over
{\delta}}\sum\limits_{j=1}^{l_{\varepsilon}}\int\limits_{0}^{\psi_{j}(\varepsilon)}\beta_{j}(\varepsilon ,x)dx \Bigr\vert <
\mbox{const.}l_{\varepsilon}\psi_{1}^{{1\over 2}}(\varepsilon)) \eqno(73)$$

\noindent is found. Let us take $\>f(\varepsilon)=\psi_{1}(\varepsilon)[\ln\psi_{1}(\varepsilon)]^{-1}$. By using the function $\>p(x)\>$ which satisfies
the condition (p1) and the inequality (42)

\begin{eqnarray*}
{1\over {\delta}}\sum\limits_{j=1}^{l_{\varepsilon}}\int\limits_{0}^{\psi_{j}(\varepsilon)}\beta_{j}(\varepsilon ,x)\ud x&>& {1\over {\delta}}
\int\limits_{0}^{\psi_{1}(\varepsilon)}\beta_{1}(\varepsilon ,x)\ud x\\
\\
&=&{1\over{3 \pi}}
\int\limits_{0}^{\psi_{1}(\varepsilon)}\sqrt{{{\alpha_{1}(x)-\varepsilon}\over{p(x)}}}\>\Big(2\alpha_{1}(x)+\varepsilon\Big)\ud x\\
\\
&>&{1\over{3\pi}}\int\limits_{{1\over{2}}f(\varepsilon)}^{f(\varepsilon)}\sqrt{{{\alpha_{1}(x)-\varepsilon}\over{p(x)}}}\>
\Big(2\alpha_{1}(x)+\varepsilon\Big)\ud x\\
\\
&>&\mbox{const.}f(\varepsilon)\Bigl(\alpha_{1}(f(\varepsilon))-\varepsilon\Bigr)^{{3\over{2}}}\qquad\qquad\qquad\qquad\quad(74)
\end{eqnarray*}

\noindent is obtained. \c{S}eng\"{u}l showed

$$\alpha_{1}(f(\varepsilon))-\varepsilon > \Bigl(\ln\psi_{1}(\varepsilon)\Bigr)^{-(\xi+1)(n+1)} \eqno(75)$$

\noindent for the small values of $\>\varepsilon >0$, \cite{Se}. From (74) and (75)

\begin{eqnarray*}
{1\over {\delta}}\sum\limits_{j=1}^{l_{\varepsilon}}\int\limits_{0}^{\psi_{j}(\varepsilon)}\beta_{j}(\varepsilon
,x)\ud x&>&\mbox{const.}{{\psi_{1}(\varepsilon)}\over{ln\psi_{1}(\varepsilon)}}\bigr(\ln\psi_{1}(\varepsilon)\bigl)^{({-3\over{2}})(\xi+1)(n+1)}\\
\\
&>&\mbox{const.}{{\psi_{1}}^{{3\over4}}}(\varepsilon) \qquad\qquad\qquad\qquad\qquad\qquad\qquad\quad(76)
\end{eqnarray*}

\noindent is found. From (73) and (76)

$$\Biggl\vert {{\sum\limits_{i=1}^{N(\varepsilon)}\lambda_{i}}\over
{{\delta}^{-1}\sum\limits_{j=1}^{l_{\varepsilon}}\int\limits_{0}^{\psi_{j}(\varepsilon)}\beta_{j}(\varepsilon ,x)dx}}-1\Biggr\vert
<\mbox{const.}l_{\varepsilon}{\psi_{1}}^{{-1\over 4}}(\varepsilon) \eqno(77)$$

\noindent is obtained. Since the series $\>\sum\limits_{m=1}^{\infty}[\alpha_{j}(0)]^{m}\>$ is convergent then we have

$$\mbox{const}>\sum\limits_{\alpha_{j}(0)\geq\varepsilon}[\alpha_{j}(0)]^{m}\geq \sum\limits_{\alpha_{j}(0)\geq\varepsilon}{\varepsilon}^{m}
={\varepsilon}^{m}l_{\varepsilon}.$$

\noindent From last inequality

$$l_{\varepsilon}<\mbox{const.}{\varepsilon}^{-m}\qquad\eqno(78)$$

\noindent is found. Since the function $\>\alpha_{1}(x)\>$ satisfy the condition {\bf $\alpha 1)$}, we have

$$\varepsilon=\alpha_{1}(\psi_{1}(\varepsilon))\geq(\ln_{n}\psi_{1}(\varepsilon))^{-\xi}\geq(\ln \psi_{1}(\varepsilon))^{-\xi}$$

\noindent for the small values of $\>\varepsilon >0$. From the last inequality above,

$$\psi_{1}(\varepsilon)>e^{{\varepsilon}^{{-1}\over{\xi}}} \qquad\eqno(79)$$

\noindent is obtained. From (77), (78) and (79)

$$\Biggl\vert {{\sum\limits_{i=1}^{N(\varepsilon)}\lambda_{i}}\over
{{\delta}^{-1}\sum\limits_{j=1}^{l_{\varepsilon}}\int\limits_{0}^{\psi_{j}(\varepsilon)}\beta_{j}(\varepsilon ,x)dx}}-1\Biggr\vert
<\mbox{const.}\varepsilon^{-m}\>e^{{{-1\over4}\varepsilon}^{-1\over {\xi}}}<\mbox{const.}e^{{-\varepsilon}^{-\beta}} \qquad\quad\eqno (80)$$

\noindent is found. We can rewrite inequality (80)

$${{\sum\limits_{i=1}^{N(\varepsilon)}\lambda_{i}}\over
{{\delta}^{-1}\sum\limits_{j=1}^{l_{\varepsilon}}\int\limits_{0}^{\psi_{j}(\varepsilon)}\beta_{j}(\varepsilon ,x)\ud x}}-
1=O\big(e^{{-\varepsilon}^{-\beta}}\big) \eqno (81)$$

\noindent as $\>\varepsilon\rightarrow 0$. From (2), (42) and (81)

$$\sum\limits_{-\lambda_{i}<-\varepsilon}\lambda_{i}={1\over {3\pi}}\biggl[
1+O\big(e^{{-\varepsilon}^{-\beta}}\big)\biggr]\sum\limits_{j}\int\limits_{\alpha_{j}(x)\geq
\varepsilon}\sqrt{{{\alpha_{j}(x)-\varepsilon}\over{p(x)}}}\>\Big(2\alpha_{j}(x)+\varepsilon\Big)\ud x$$

\noindent as $\>\varepsilon \rightarrow 0$, is obtained. $\Box$

\noindent Let us assume that the function $\>\alpha_{1}(x)\>$ satisfies the following condition:

\noindent {\bf$\alpha 2)$} For every $\>\eta >0\>$

$$\lim_{x\rightarrow\infty}\alpha_{1}(x)x^{a_{0}-\eta}=\lim_{x\rightarrow\infty}[\alpha_{1}(x)x^{a_{0}+\eta}]^{-1}=0$$

\noindent Here, $\>a_{0}\>$ is a constant in the interval $\>(0,\frac{2}{3})$.

\begin{theorem}
We suppose that the operator function $\>Q(x)\>$, the scalar function $\>p(x)\>$ satisfy the condition $\>Q1)-Q3),\>p1)-p3)\> $ and
$\>\alpha_{1}(x)\>$ also satisfies the condition {\bf$\alpha 2)$} . In addition the series $\>\sum\limits_{j=1}^{\infty}[\alpha_{j}(0)]^{m}\>$ is
convergent for a constant $\>m\>$ satisfying the condition

$$0< m <{{(2-3a_{0})^{2}}\over{2a_{0}(4-3a_{0})}} \qquad \eqno(82)$$

\noindent then the asymptotic formula

$$ \sum\limits_{i=1}^{N(\varepsilon)} \lambda_{i}=\frac{1}{3\pi}[1+O({\varepsilon}^{t_{0}})] \sum\limits_{j}\int\limits_{\alpha_{j}(x)\geq\varepsilon}
\sqrt{{{\alpha_{j}(x)-\varepsilon}\over{p(x)}}}\>\Big(2\alpha_{j}(x)+\varepsilon\Big)dx $$

\noindent is satisfied as $\>\varepsilon\rightarrow 0$. Where $\>t_{0}\>$ is a positive constant.
\end{theorem}

\noindent {\bf Proof :} By Theorem 2.4 and Theorem 2.5, we have

$$\Biggl\vert \sum\limits_{i=1}^{N(\varepsilon)}\lambda_{i}-{\delta}^{-1}\sum\limits_{j=1}^{l_{\varepsilon}}
\int\limits_{0}^{\psi_{j}(\varepsilon)}\beta_{j}(\varepsilon ,x)dx \Biggr\vert < \mbox{const.}l_{\varepsilon}
\Bigl(\int\limits_{0}^{\delta}\alpha_{1}^{\frac{3}{2}}(x)dx +\psi_{1}^{a}(\varepsilon)\Bigr)\eqno(83)$$

\noindent for the small values of $\>\varepsilon >0$. Since the function $\>\alpha_{1}(x)\>$ is decreasing,

$$ \alpha_{1}(x)\geq\alpha_{1}(\psi_{1}(2\varepsilon))=2\varepsilon \eqno(84)$$

\noindent in the interval $\>[0,\psi_{1}(2\varepsilon)]$. Since the function $\>p(x)\>$ satisfies the condition  p1) and (42), (84) then we find

\begin{eqnarray*}
\delta^{-1}\sum\limits_{j=1}^{l_{\varepsilon}}\int\limits_{0}^{\psi_{j}(\varepsilon)}\beta_{j}(\varepsilon ,x)\ud x
&>&\frac{1}{3\pi}\int\limits_{0}^{\psi_{1}(\varepsilon)} \sqrt{{{\alpha_{1}(x)-\varepsilon}\over{p(x)}}}\>\Big(2\alpha_{1}(x) +\varepsilon\Big)\ud x\\
\\
&>&\mbox{const.}{\varepsilon}^{\frac{3}{2}}\psi_{1}(2\varepsilon)  \qquad\qquad\qquad\qquad\qquad\qquad\quad(85)
\end{eqnarray*}

\noindent If we consider that the function $\>\alpha_{1}(x)\>$  satisfies the condition {\bf$\alpha 2)$} and

$\>\lim_{\varepsilon \rightarrow 0} \psi_{1}(\varepsilon)=\infty$, then we have

$$\lim_{\varepsilon\rightarrow\infty}[\alpha_{1}(\psi_{1}(2\varepsilon))(\psi_{1}(2\varepsilon))^{a_{0}+\eta}]^{-1}=0$$

\noindent From the last equality above, we obtain

$$\psi_{1}(2\varepsilon)>(\varepsilon)^{\frac{-1}{a_{0}+\eta}} \eqno(86)$$

\noindent for the small value of $\>\varepsilon >0$. From (85) and (86)

$$ {\delta}^{-1}\sum\limits_{j=1}^{l_{\varepsilon}}\int\limits_{0}^{\psi_{j}(\varepsilon)}\beta_{j}(\varepsilon ,x)\ud x
>\mbox{const.}{\varepsilon}^{\frac{3a_{0}+3\eta-2}{2(a_{0}+\eta)}}  \eqno(87)$$

\noindent is found. We limit the integral $\>\int\limits_{0}^{\delta}\alpha_{1}^{\frac{3}{2}}(x)\ud x \>$ at the right hand side of the inequality
(83). Since the function $\>\alpha_{1}(x)\>$  satisfies the condition {\bf$\alpha 2)$}, then we have
$$\>\alpha_{1}(x)\leq\mbox{const.}x^{\eta-a_{0}} \qquad (\eta <a_{0}). \eqno(88)$$

\noindent Therefore we have

$$\int\limits_{0}^{\delta}\alpha_{1}^{\frac{3}{2}}(x)dx \leq\mbox{const.}\int\limits_{0}^{\delta}x^{{\frac{3}{2}}(\eta -a_{0})}\ud x
<\mbox{const.}\delta^{{\frac{1}{2}}(2-3a_{0}+3\eta)}.\eqno(89)$$

\noindent On the other hand, from (3)

$$\delta< \psi_{1}^{1-a}(\varepsilon) \eqno(90)$$

\noindent is obtained. If we take $\>x=\psi_{1}(\varepsilon)\>$ in the inequality (88), then we find

$$\>\alpha_{1}(\psi_{1}(\varepsilon))\leq\mbox{const.}\psi_{1}^{\eta-a_{0}}(\varepsilon) \qquad (\eta <a_{0}) $$

\noindent or
$$\psi_{1}(\varepsilon)\leq\mbox{const.}{\varepsilon}^{\frac{-1}{a_{0}-\eta}} \eqno(91)$$

\noindent From (89), (90) and (91), we have

$$\int\limits_{0}^{\delta}\alpha_{1}^{\frac{3}{2}}(x)dx \leq\mbox{const.}{\varepsilon}^{-\frac{(1-a)(2-3a_{0}+3\eta)}{2(a_{0}-\eta)}}.
\eqno(92)$$

\noindent From (78), (91) and (92)

$$l_{\varepsilon}\int\limits_{0}^{\delta}\alpha_{1}^{\frac{3}{2}}(x)\ud x<
\mbox{const.}{\varepsilon}^{-m-\frac{(1-a)(2-3a_{0}+3\eta)}{2(a_{0}-\eta)}} \eqno(93)$$

$$l_{\varepsilon}\psi_{1}^{a}(\varepsilon)<\mbox{const.}{\varepsilon}^{\frac{-m(a_{0}-\eta)+a}{(a_{0}-\eta)}} \eqno(94)$$

\noindent are found.
\noindent From (87), (93) and (94) we obtain

$${{l_{\varepsilon}\int\limits_{0}^{\delta}\alpha_{1}^{\frac{3}{2}}(x)\ud x}\over
{\delta^{-1}\sum\limits_{j=1}^{l_{\varepsilon}}\int\limits_{0}^{\psi_{j}(\varepsilon)}\beta_{j}(\varepsilon ,x)\ud x}}
< \mbox{const.}{\varepsilon}^{F_{1}(\eta)} \eqno(95)$$
\noindent and
$${{l_{\varepsilon}\psi_{1}^{a}(\varepsilon)\ud x}\over
{\delta^{-1}\sum\limits_{j=1}^{l_{\varepsilon}}\int\limits_{0}^{\psi_{j}(\varepsilon)}\beta_{j}(\varepsilon ,x)\ud x}}
< \mbox{const.}{\varepsilon}^{F_{2}(\eta)}. \eqno(96)$$

\noindent Here,

$$ F_{1}(\eta)=-m-{{(1-a)(2-3a_{0}+3\eta)}\over{2(a_{0}-\eta)}}-{{3a_{0}+3\eta-2}\over{2(a_{0}+\eta)}}$$

$$ F_{2}(\eta)=-{{m(a_{0}-\eta)+a}\over{(a_{0}-\eta)}}-{{3a_{0}+3\eta-2}\over{2(a_{0}+\eta)}}.$$

\noindent There is a number $\>\omega=\omega(t)>0 \qquad(0<\eta < \omega)\>$ such that

$$ F_{1}(\eta)>{{2a-2a_{0}m-3aa_{0}}\over{2a_{0}}}-t \eqno(97)$$

$$ F_{2}(\eta)>{{2-3a_{0}-2a_{0}m-2a}\over{2a_{0}}}-t \eqno(98)$$

\noindent for every $\>t>0\>$. If we take

$$a={{(2-3a_{0})^{2}+6a_{0}^{2}m}\over{4(2-3a_{0})}}, \qquad t=t_{0}={1\over{16a_{0}}}\Bigl((2-3a_{0})^{2}+6a_{0}^{2}m-8a_{0}m\Bigr)$$

\noindent in the inequalities (97) and (98), then we have

$$F_{1}(\eta)>t_{0}\>\quad;\>\quad  F_{2}(\eta)>t_{0}. \eqno(99)$$

\noindent Since the number $\>m\>$ satisfies the condition (82), we have $\>a\in(0,1)\>$ and $\>t_{0}>0$.

\noindent From (83), (95),(96) and  (99) we obtain

$$\Biggl\vert {{\sum\limits_{i=1}^{N(\varepsilon)}\lambda_{i}}\over
{{\delta}^{-1}\sum\limits_{j=1}^{l_{\varepsilon}}\int\limits_{0}^{\psi_{j}(\varepsilon)}\beta_{j}(\varepsilon ,x)dx}}-1\Biggr\vert
<\mbox{const.}{\varepsilon}^{t_{0}} .\eqno (100)$$

\noindent By (42), (97) and (100) we have the asymptotic formula

$$ N(\varepsilon)=\frac{1}{3\pi}\Bigl[1+O({\varepsilon}^{t_{0}})\Bigr]\> \sum\limits_{j}\int\limits_{\alpha_{j}(x)\geq\varepsilon}
\sqrt{{{\alpha_{j}(x)-\varepsilon}\over{p(x)}}}\>\Bigl(2\alpha_{j}(x)+\varepsilon\Bigr)\ud x $$

\noindent as $\>\varepsilon\rightarrow 0$. $\Box$

\begin{exam}
Let $\>H =L_{2}[0,\pi ]\>$ be a separable Hilbert space and \\ $\>e_{i}=\sqrt{\frac{2}{\pi}}\>\sin ix \quad (i=1,2,\cdots)\>$ be a standard basis in $\>H\>$. Let $\>Q(x):H \rightarrow H $
$$ Q(x)f=\sum \limits_{i=1}^{\infty} \alpha(x)i^{-2}\Big(f,e_{i}\Big)e_{i} \quad (f\in H)$$

\noindent for all $\>x \in [0,\infty)$. $\>Q(x)\>$ is a self adjoint, completely continuous and positive operator function.
\noindent The eigenvalues of $\>Q(x)\>$ are in the form
\begin{displaymath}
\alpha (x)=\left\{\begin{array}{ll}
\frac{2}{\ln\ln b}-\frac{x}{b\ln\ln b}\>,  &0 \leq x \leq b\\
\\
\frac{1}{\ln\ln x}\>, &b \leq x < \infty
\end{array} \right.
\end{displaymath}

\noindent  Here $\>b>e^{3}\>$ is a constant such that $\>\ln x > (\ln\ln x)^{2}$.
\end{exam}


\begin{thebibliography}{13}

\bibitem{Se} \c{S}eng\"{u}l, S. \emph{The asymptotic behaviour of the spectrum of negative
part of Sturm-Liouville problem with operator coefficient }, PhD thesis in YT\"{U} FBE (2006).
(In Turkish).
$$https://tez.yok.gov.tr/UlusalTezMerkezi/TezGoster?key=$$
$$-L8ilcwn9ZRRc_YMKxXW1u5yjU0aBL-ngSBJoTERxJektuxRXLjFdBQGCdCgtfhR$$

\bibitem{Ad} Ad{\i}g\"{u}zelov, E.E. , Oer , Z. \emph{Asymptotic Expansion for
the sum of negative Eigenvalues of Sturm-Lioville operator given in Semi-axis},
YT\"UD, (2002), Vol 1,26-35.

\bibitem{Ba} Bak\c{s}i, \"{O}., Ismay{\i}lov,S. \emph{An asymptotic formula for the sum of
negative eigenvalues of second order differntial operator given in infinite interval},
Sigma M\"{u}hendislik ve Fen Bilimleri Dergisi
2005-4, 87-98. (In Turkish)

\bibitem{Sk}Ska\d{c}ek B.Y. \emph{ Asymptod of Negative Part of Spectrum of One Dimensioned
Differential Operators}, Pribl. metodi res.eniya differn. uraveniy, Kiev, 1963",
Pribl. Metod reseniya differens, unavneniy, Kiev, (1963).

\bibitem{Ae} Ad{\i}g\"{u}zelov, E.E. \emph{The asymptotic behaviour of the spectrum's
negative part of Sturm-Liouville problem with operator coefficient}, Izv. AN Az.SSR,
Seriya fiz.-tekn.i mat. nauk, No:6, 8-12, (1980). (In Russian)

\bibitem{Ma} Maksudov F.G.,
Bayramo\v{g}lu M., Ad{\i}g\"{u}zelov E.,\emph{On asymptotics of spectrum and trace
of high order differantial operator with operator coefficients},
Do\v{g}a-Turkish journal of Mathematics, (1993), vol.17.

\bibitem{Ab} Ad{\i}g\"{u}zelov, E.E., Bak\c{s}i, \"{O}., Bayramov, A.M.
\emph{The Asymptotic Behaviour of the Negative Part of the Spectrum of
Sturm-Liouville Operator with the Operator Coefficient which Has Singularity},
International Journal of Differential Equations and Applications,
Vol.6, No.3, 315-329, (2002).

\bibitem{Co} Gohberg, I.C. and Krein, M.G., \emph{Introduction to the Theory of Linear
Non-self Adjoint Operators in Hilbert Space}, Translation of Mathematical Monographs,
Vol.18 (AMS, Providence, R.I.,1969).

\bibitem{Ly} Lysternik, L.A. and Sobolev, V.I. \emph{Elements of Functional Analysis},
(English translation ), John Willey Sons, New York, page 229, (1974).

\bibitem{Gl} Glazman, I.M. \emph{Direct methods qualitative spectral analysis of singular
differntial operators}, Jerusalem, pages 34-44, (1965).

\bibitem{Sm} Smirnov,V.I., \emph{A Course of Higher Mathematics}, vol.5, Pergamon Pres, New York, page 623, (1964).

\bibitem{Ch} Courant, R. and Hilbert, D., \emph{Methods of Mathematical Physics},
vol.1, New York, page 408, (1966).


\end{thebibliography}
 \end{document}